\documentclass{article}
\usepackage{arxiv}
\usepackage[english]{babel}
\usepackage[T1]{fontenc}
\usepackage[utf8]{inputenc}
\usepackage{graphicx}
\usepackage{amssymb}
\usepackage[numbers]{natbib}
\usepackage{hyperref}
\usepackage[ruled,vlined]{algorithm2e}

\makeatletter
\def\blfootnote{\xdef\@thefnmark{}\@footnotetext}
\makeatother

\begin{document}
\title{Krylov Iterative Methods for the Geometric Mean of Two Matrices Times a Vector}
\author{Jacopo Castellini\\
Department of Mathematics and Computer Science\\
Università degli Studi di Perugia\\
\texttt{jacopo.castellini@studenti.unipg.it}}
\date{}
\maketitle

\begin{abstract}
In this work, we are presenting an efficient way to compute the geometric mean of two positive definite matrices times a vector. For this purpose, we are inspecting the application of methods based on Krylov spaces to compute the square root of a matrix. These methods, using only matrix-vector products, are capable of producing a good approximation of the result with a small computational cost.
\blfootnote{This work is published in \textit{Numerical Algorithms 74(2), 561--571, Springer, 2017}}
\end{abstract}


\section{Introduction}
Many problems in engineering, medicine and computer science make use of the geometric mean of two matrices $A\#B=A(A^{-1}B)^{1/2}$, with $A$ and $B$ being positive definite matrices. These applications include the calculation of electrical networks \cite{trapp}, diffusion tensor imaging \cite{moakher} and image deblurring \cite{benedetto}. For these applications there are methods to approximate a function of a matrix, as shown in \cite{iannazzo} and \cite{triplenick}, that are often used with acceptable results in term of computational cost.

But in the solution of elliptic partial differential equations with domain decomposition (such as in \cite{arioli} and \cite{kourounis}), in certain cases it is required to compute the geometric mean of two matrices times a vector $(A\#B)v$, where $A$ and $B$ are often very large and sparse. A matrix of dimension $n\times m$ is called sparse if defined $\mathcal{S}(n,m)$ as the number of its non-zero elements, it holds
\begin{equation}
\lim_{n,m\rightarrow \infty}\frac{\mathcal{S}(n,m)}{nm}=0.
\end{equation}

The use of algorithms created to compute $A\#B$ could lead to unsustainable computational times or sometimes lead to the impossibility of calculating the result, since these methods do not exploit the sparse structure of data, that instead should represent an advantage and a saving in terms of operations.

For these reasons, methods based on Krylov spaces are preferred. In \cite{arioli} and \cite{parlett}, the usage of a generalized version of the Lanczos method is exploited (we briefly present it in Section 3). In this work, we will use methods based on both polynomial and rational Krylov spaces in order to compute the geometric mean of two matrices times a vector in an efficient way, by exploiting the sparse structure of the involved matrices. In fact, these methods are easier to use, because they only compute matrix-vector products and solve linear systems with the conjugate gradient method. Moreover, their use as iterative methods is another advantage to reduce the computational cost. Specifically, we will focus on the Arnoldi method and on some variations of the rational Arnoldi method, explained in Section \ref{polyarnoldi} and Section \ref{ratarnoldi} respectively. Experiments presented in Section \ref{experiments} show how these methods, especially rational Arnoldi ones, are efficient in terms of approximation and computational time.

\section{Geometric Mean of Two Matrices}
The geometric mean $A\#B$ of two positive definite matrices $A,B\in\mathbb C^{n\times n}$ is defined in \cite{bhatia} as
\begin{equation}
A(A^{-1}B)^{1/2}.
\end{equation}

Let $M$ be a diagonalizable matrix with positive eigenvalues, so that an invertible matrix $K$ and a diagonal matrix $D=\mbox{diag}({\lambda}_1,\ldots,{\lambda}_n)$ exist such that $M=KDK^{-1}$. The square root of M is defined as
\begin{equation}
M^{1/2}=K\mbox{diag}(\sqrt{{\lambda}_1},\ldots,\sqrt{{\lambda}_n})K^{-1}
\end{equation}
that turns out to be still positive definite. The matrix $A^{-1}B$ is diagonalizable with positive eigenvalues, because it is similar to $(A^{1/2})^{-1}B(A^{1/2})^{-1}$ which is positive definite, so the notation $A(A^{-1}B)^{1/2}$ makes sense. It can be shown that $(A^{1/2})^{-1}=(A^{-1})^{1/2}$, then we use the notation $A^{-1/2}$ without ambiguity. We observe also that
\begin{equation}
A\#B=A(A^{-1}B)^{1/2}=(BA^{-1})^{1/2}A=B(B^{-1}A)^{1/2}B=B\#A.
\end{equation}

We also have that, if both $A$ and $B$ are positive definite matrices, then so is $A(A^{-1}B)^{1/2}$.

To get an efficient computation of $(A\#B)v$ it is important not to use expensive operations. The most problematic part is to find a good approximation of $A(A^{-1}B)^{1/2}v=A(B^{-1}A)^{-1/2}v$ without explicitly computing $B^{-1}$ or the square root of $B^{-1}A$. We recall the fact that $z^{-1/2}$ is a Markov function, i.e. a function of the form
\begin{equation}
f(z)=\int_{\Gamma}\frac{d\gamma(x)}{z-x},
\end{equation}
with $\gamma$ complex measure over the closed set $\Gamma\subseteq\mathbb{C}$, because
\begin{equation}
f(z)=z^{-1/2}=\int_{-\infty}^{0}\frac{1}{z-x}\frac{dx}{\pi\sqrt{-x}}.
\end{equation}

\section{Polynomial Krylov Spaces}
\label{polyarnoldi}
Given $A\in\mathbb{C}^{n\times n}$ and $b\in\mathbb{C}^n$, the $m$-th Krylov polynomial space $K_m(A,b)$ associated with them is defined as
\begin{equation}
K_m(A,b)=\mbox{span}\{b,Ab,A^2b,\ldots ,A^{m-1}b\}\subseteq\mathbb{C}^n.
\end{equation}

With the increase of $m$, we get Krylov spaces all nested one inside another, $K_1(A,b)\subseteq K_2(A,b)\subseteq\ldots\subseteq K_d(A,b)$. It can be shown that if $b\ne 0$, there exists $k\le n$ such that
\begin{equation}
\mbox{dim }K_i(A,b)=\left\{
\begin{array}{ll}
i & \mbox{if }i<k, \\
k & \mbox{otherwise}.
\end{array}
\right.
\end{equation}

For each step the space $K_i(A,b)$ has a basis $\{b,Ab,\ldots,A^{i-1}b\}$, but from a certain $k+1$ on the new term $A^kb$ is linearly dependent with the other vectors $b,\ldots,A^{k-1}b$ in the basis, giving that $K_{k+1}(A,b)$ is the same of $K_k(A,b)$.

\subsection{Generalized Lanczos Method}
First we briefly recall the variation of the Lanczos method presented in \cite{parlett} and \cite{arioli}, called the generalized Lanczos method (see \cite{parlett}, \cite{bini} and \cite{golub} for an introduction to the Lanczos method). This iterative method is adapted in \cite{arioli} in order to compute the geometric mean of two positive definite matrices times a vector. Given two positive definite matrices $A,B\in\mathbb{C}^{n\times n}$, the method constructs at the step $k$ a new vector $q_k\in\mathbb{C}^n$ that is $B$-orthogonal, that is $q_i^*Bq_j=0$ with $i\ne j$, such that for each step we have
\begin{equation}
AQ_k=BQ_kT_k+{\beta}_{k+1}Bq_{k+1}e_k^T\qquad Q_k^*BQ_k=I_k,
\end{equation}
with $Q_k=[q_1,\ldots,q_k]$, $I_k$ the $k\times k$ identity matrix and $e_k$ its $k$-th column.

The matrix $T_k\in\mathbb{C}^{k\times k}$ is the projection of $A$ into the space generated by the columns of the matrix $Q_k$. At the $n$-th step this gives us the equalities
\begin{equation}
Q^*AQ=T,\qquad Q^*BQ=I.
\end{equation}

\begin{algorithm}[htbp]
\TitleOfAlgo{Generalized Lanczos Method}
${\beta}_0=0$\;
$q_1=v/\|v\|_B$\;
\For{$i=1\ldots k$}{
$w=B^{-1}Aq_i-{\beta}_{i-1}q_{i-1}$\;
${\alpha}_i=w^*Bq_i$\;
$w=w-{\alpha}_iq_i$\;
${\beta}_i=\|v\|_B$\;
$q_{i+1}=w/{\beta}_i$\;
}
$Q=[q_1,\ldots,q_k]$\;
$T=\mbox{tridiag}(\beta,\alpha,\beta)$\;
$(A\#B)v=BQT^{1/2}e_1\|v\|_B$\;
\end{algorithm}

Using this method we find an expression for $B\#A$ (that is the same of $A\#B$) as
\begin{equation}
B\#A=B(B^{-1}A)^{1/2}=(QT^{-1/2}Q^*)^{-1}=BQT^{1/2}Q^*B.
\end{equation}

To compute $(B\#A)v$, with $v\in\mathbb{R}^n$, we can stop at the $k$-th step, and taking $q_1=v/\|v\|_B$, we have
\begin{equation}
(B\#A)v\approx BQ_kT_k^{1/2}e_1\|v\|_B,
\end{equation}
because $Q_k^*Bv=e_1\|v\|_B$, that is a better approximation the closer $k$ is to $n$.

\subsection{Arnoldi Method}
The Arnoldi method, exposed in \cite{higham}, is a general version of the Lanczos method that does not need the starting matrices to be positive definite. This method computes the decomposition
\begin{equation}
\label{arnoldi}
Q^*AQ=H,
\end{equation}
where $Q\in\mathbb{C}^{n\times n}$ is unitary, and the vectors $q_1,\ldots,q_n$ are its columns, and $H\in\mathbb{C}^{n\times n}$ is upper Hessenberg. From the equality (\ref{arnoldi}) we get that
\begin{equation}
Aq_k=\sum_{i=1}^{k+1}h_{ik}q_i,\qquad\mbox{with }k=1,\ldots,n-1.
\end{equation}

The previous relation can be rewritten as
\begin{equation}
h_{k+1,k}q_{k+1}=Aq_k-\sum_{i=1}^kh_{ik}q_i=r_k,
\end{equation}
and, because $Q$ is unitary, we have that
\begin{equation}
h_{ik}=q_i^*Aq_k\qquad\mbox{with }i=1,\ldots,k.
\end{equation}

If $r_k\ne 0$, so $q_{k+1}=r_k/h_{k+1,k}$, where $h_{k+1,k}=\|r_k\|_2$. The Arnoldi method is a Krylov space method because
\begin{equation}
\mbox{span}\{q_1,q_2,\ldots,q_k\}=\mbox{span}\{q_1,Aq_1,\ldots,A^{k-1}q_1\},
\end{equation}
that is the vectors $q_1,\ldots,q_k$ constructed by the Arnoldi method are a basis for the Krylov space $K_k(A,b)$. Used as an iterative method, at the $k$-th step the Arnoldi method give us the factorization
\begin{equation}
AQ_k=Q_kH_k+h_{k+1,k}q_{k+1}e_k^T,
\end{equation}
where $Q_k$ consist of the first $k$ column of $Q$, $H_k$ is a $k\times k$ upper Hessenberg matrix and $e_k$ is the $k$-th column of the $k\times k$ identity matrix. Because
\begin{equation}
Q_k^*AQ_k=H_k,
\end{equation}
we have that $H_k$ is the projection of $A$ in $K_k(A,b)$. The orthonormalization of vectors $q_1,\ldots,q_n$ comes from the modified Gram-Schmidt method (see \cite{trefethen}), but in finite precision arithmetic a loss of orthogonality can occur, as shown in \cite{parlett}.

Using the Arnoldi method we can approximate $f(A)b$, starting with the vector $q_1=b/\|b\|_2$. As shown in \cite{higham} and \cite{faber} it stands that, at the $k$-th step, the approximation is
\begin{equation}
f_k=Q_kf(H_k)e_1\|b\|_2=Q_kf(H_k)Q_k^*b,
\end{equation}
that is the same of calculating $f$ into the Krylov space $K_k(A,B)$ and then expanding the result into the original space $\mathbb{C}^{n\times n}$. So an approximation for $(A\#B)v=A(A^{-1}B)^{1/2}v$ is
\begin{equation}
(A\#B)v\approx AQ_kH_k^{-1/2}Q_k^*v=AQ_kH_k^{-1/2}e_1\|v\|_2.
\end{equation}

\begin{algorithm}[H]
\TitleOfAlgo{Arnoldi Method}
$q_1=v/\|v\|_2$\;
\For{$i=1\ldots k$}{
$w=(A^{-1}B)q_i$\;
\For{$j=1\ldots i$}{
$h_{j,i}=q_j^*w$\;
$z=z-h_{j,i}q_j$\;
}
$h_{i+1,i}=\|z\|_2$\;
$q_{i+1}=w/h_{i+1,i}$\;
}
$Q=[q_1,\ldots,q_k]$\;
$H=\{h_{i,j}\}_{i,j=1,\ldots,k}$\;
$(A\#B)v=AQH^{1/2}Q^*v$\;
\end{algorithm}

\section{Rational Krylov Spaces}
\label{ratarnoldi}
The definition of rational Krylov spaces is similar to the polynomial one, except for the presence of a denominator (see \cite{guttel1} and \cite{guttel2} for an introduction). Given $A\in\mathbb{C}^{n\times n}$ and $b\in\mathbb{C}^n$, and the sequence of polynomials $q_{m-1}(A)$ of degree $m-1$
\begin{equation}
q_{m-1}(z)=\prod_{j=1}^{m-1}(1-z/{\xi}_j),
\end{equation}
in which the values ${\xi}_j\in\mathbb{C}\cup\{\infty\}$ are called poles and are numbers in the extended complex plane different from all the eigenvalues of $A$ and $0$, the rational Krylov space of order $m$ associated to them is defined as
\begin{equation}
\label{denominator}
Q_m(A,b)=q_{m-1}(A)^{-1}\mbox{span}\{b,Ab,\ldots,A^{m-1}b\}\subseteq\mathbb{C}^n.
\end{equation}

In the case in which ${\xi}_j=\infty$ for some $j$ the corresponding factor $(1-z/{\xi}_j)$ in (\ref{denominator}) is replaced by $1$.

In the previous definition it is not taken into account the case ${\xi}_j=0$, but we can exclude every other value $\sigma$ simply by using the new values $\widehat{A}=A-\sigma I$ and $\widehat{{\xi}_j}={\xi}_j-\sigma$. This can be done considering the new polynomial $q_{m-1}$ that uses poles $\widehat{\xi}_j$ and the new matrix $\widehat{A}$ instead of $A$. Like the polynomial ones, rational Krylov spaces are also of increasing dimension and are nested until a certain dimension $k$ from which they do not change anymore. Polynomial Krylov spaces can be seen as a special case of the rational ones, when every pole ${\xi}_j$ is $\infty$, that is when $q_{m-1}=1$.

\subsection{Rational Arnoldi Method}
To compute a basis $V_m=[v_1,\ldots ,v_m]\in\mathbb{C}^{n\times m}$ for the rational Krylov space $Q_m(A,b)$ we can use the rational Arnoldi method, as proposed in \cite{guttel1} and in \cite{faber}. Starting with $v_1=b/\|b\|_2$, in the following iterations the vector $v_{j+1}$ is generated orthogonalizing
\begin{equation}
x_j=(I-A/{\xi}_j)^{-1}Av_j
\end{equation}
against the previous orthonormal vectors $v_1,\ldots ,v_j$. If $h_{j+1,j}\ne 0$, we have that
\begin{equation}
x_j=\sum_{i=1}^{j+1}v_ih_{i,j},
\end{equation}
and $v_{j+1}=x_j/h_{j+1,j}$, instead if $h_{j+1,j}=0$ we can proceed finding a new vector $v_{j+1}$ that is orthonormal to $v_1,\ldots,v_n$. 

We also have that
\begin{equation}
A\left(v_j+\sum_{i=1}^{j+1}v_ih_{i,j}{\xi}_j^{-1}\right)=\sum_{i=1}^{j+1}v_ih_{i,j},
\end{equation}
that gives us the decomposition
\begin{equation}
AV_m(I_m+H_mD_m)+Av_{m+1}h_{m+1,m}{\xi}_m^{-1}e_m^T=V_mH_m+v_{m+1}{\xi}_m^{-1}e_m^T,
\end{equation}
in which $D_m=\mbox{diag}({\xi}_1^{-1},\ldots,{\xi}_m^{-1})$, $I_m$ is the $m\times m$ identity matrix and $e_m$ is its $m$-th column. Defining
\begin{equation}
\underline{H_m}=\left[
\begin{array}{c}
H_m \\
h_{m+1,m}e_m^T
\end{array}
\right]\qquad\mbox{and}\qquad\underline{K_m}=\left[
\begin{array}{c}
I_m+H_mD_m \\
h_{m+1,m}{\xi}_m^{-1}e_m^T
\end{array}
\right]
\end{equation}
we have that
\begin{equation}
AV_{m+1}\underline{K_m}=V_{m+1}\underline{H_m},
\end{equation}
where $V_{m+1}=[V_m,v_{m+1}]$. Finally, if the last pole ${\xi}_m$ is $\infty$, the decomposition is simplified into
\begin{equation}
AV_mK_m=V_{m+1}\underline{H_m},
\end{equation}
where $K_m$ is the $m\times m$ upper part of $\underline{K_m}$.

Computed the basis $V_m$ of $Q_m(A,b)$, we can approximate $f(A)b$ as
\begin{equation}
f_m^{RA}=V_mf(A_m)V_m^*b,\qquad\mbox{where } A_m=V_m^*AV_m\in\mathbb{C}^{m\times m}.
\end{equation}

\begin{algorithm}[htbp]
\TitleOfAlgo{Rational Arnoldi Method}
$q_1=v/\|v\|_2$\;
\For{$i=1\ldots k$}{
$w=(I-B^{-1}A/{\xi}_i)^{-1}B^{-1}Aq_i$\;
\For{$j=1\ldots i$}{
$h_{j,i}=q_j^*w$\;
$z=z-h_{j,i}q_j$\;
}
$h_{i+1,i}=\|z\|_2$\;
$q_{i+1}=w/h_{i+1,i}$\;
}
$V=[q_1,\ldots,q_k]\;$
$H=\{h_{i,j}\}_{i,j=1,\ldots,k}$\;
$K=I+H\mbox{ diag}({\xi}_1^{-1},\ldots,{\xi}_k^{-1})$\;
$(A\#B)v=AV(HK^{-1})^{-1/2}V^*v$\;
\end{algorithm}

The benefit of this method is that, in many interesting cases, $f_m^{RA}$ is a good approximation of $f(A)b$ also if $m$ is small: it is required to compute $f(A_m)$, but $A_m$ is small compared to $A$ itself. If the last pole ${\xi}_m$ is $\infty$, it is not even necessary to compute $A_m=V_m^*AV_m$, but it can be computed $A_m=H_mK_m^{-1}$. Also $H_m$ and $K_m$ are small compared to $A$, so computing $K_m^{-1}$ or the product $H_mK_m^{-1}$ is possible. For our problem $(A\#B)v$ we have found the approximation
\begin{equation}
(A\#B)v\approx AV_m(H_mK_m^{-1})^{-1/2}V_m^*v=AV_m(H_mK_m^{-1})^{-1/2}e_1\|v\|_2,
\end{equation}
starting the method with $v_1=v/\|v\|_2$.

From a Crouzeix theorem exposed in \cite{crouzeix} we know that there exists a constant $C\leq 11.08$ for which $\|f(A)\|_2\leq C\|f\|_{\Sigma}$, where the second norm indicates the maximum absolute value of $f$ over the compact set $\Sigma$ on which the function is approximated, that in this case is a set containing the spectrum of $A$. As explained in \cite{guttel1}, this helps us to find a result of semi-optimality for the approximation given by the rational Arnoldi method: in fact, if $f$ is analytic in $\Sigma$ and we define $f_m=r_m(A)b$ with $r_m\in\mathcal{P}_{m-1}/q_{m-1}$, where $\mathcal{P}_{m-1}$ is a polynomial of degree $m-1$, so it stands
\begin{equation}
\|f(A)b-f_m^{RA}\|_2\leq 2C\|b\|_2\min_{r_m\in\mathcal{P}_{m-1}/q_{m-1}}\|f-r_m\|_{\Sigma}.
\end{equation}

It seems obvious that the choice of the poles ${\xi}_j$ is strongly connected to the function $f$ which we have to approximate. To find a small value of $\|f(A)b-f_m^{RA}\|_2$, it is necessary to search a sufficiently uniform approximation of $f$ over $\Sigma$. Moreover, sometimes it is convenient to consider only a certain set $\Xi$ for the poles disjoint from $\Sigma$. For Markov functions, and especially for the $z^{-1/2}$ function that appears in our problem, a good set $\Xi$ is $[-\infty,0]$, and as $\Sigma$ we can use $[{\lambda}_{min},{\lambda}_{max}]$, with ${\lambda}_{min}$ and ${\lambda}_{max}$ respectively the minimum and maximum eigenvalues of $A$ (as above, see \cite{guttel1} for further informations). Computing the eigenvalues of $A$ can be computationally expensive, however only the extreme eigenvalues are required.

\subsubsection{Extended Krylov Method}
A special case of the rational Arnoldi method is the so-called extended Krylov method, proposed in \cite{druskin} and then investigated in \cite{simoncini1} and \cite{simoncini2}. It simply chooses alternately the poles as ${\xi}_{even}=\infty$ and ${\xi}_{odd}=0$. This method is good for approximate Markov functions. Its benefit is obvious: the choice of the poles is completely a priori, without information on the spectrum of $A$. It had been observed in \cite{druskin} that in some cases this is equivalent to a rational Arnoldi method in which all the poles are chosen as an unique asymptotically optimal value. However, finding the optimal pole needs some informations on the spectrum of $A$.

\subsubsection{Generalized Leja Points}
An alternative way to choose poles, exposed in \cite{guttel1}, is the so-called generalized Leja points method, that rely on some logarithmic potential theory's instruments (for a complete mathematical explanation about it see \cite{ransford} and \cite{saff}). Given two closed sets $\Sigma$ and $\Xi$, both of non-zero logarithmic capacity (the logarithmic capacity of a set is defined as $\mbox{cap}(E)=e^{V_E}$, where $V_E$ is the Robin constant for $E$) and of positive distance, the pair $(\Sigma ,\Xi)$ is called a condenser, and we can associate a number to it called condenser capacity $\mbox{cap}(\Sigma ,\Xi)$ (still see \cite{saff} for further explanations). We now consider the sequence of functions
\begin{equation}
s_m(z)=\frac{(z-{\sigma}_1)\cdots(z-{\sigma}_m)}{(1-z/{\xi}_1)\cdots(1-z/{\xi}_m)},\qquad\mbox{with }m=1,\ldots,n,
\end{equation}
with nodes ${\sigma}_j\in\Sigma$ and poles ${\xi}_j\in\Xi$. Our aim is to make this sequence in absolute value as large as possible over $\Xi$ and as small as possible over $\Sigma$. It can be shown that for this kind of sequences the relation
\begin{equation}
\lim_{m\rightarrow\infty}\sup\left(\frac{{\sup}_{z\in\Sigma}|s_m(z)|}{{\inf}_{z\in\Xi}|s_m(z)|}\right)^{1/m}\geq e^{-1/\mbox{cap}(\Sigma,\Xi)}.
\end{equation}
lasts.

To find sequences of functions in which this inequality is true is called the generalized Zolotarev problem for the condenser $(\Sigma ,\Xi)$.

A practical method to obtain such functions is the following greedy algorithm: starting with ${\sigma}_1$ and ${\xi}_1$ as points of minimum distance from $\Sigma$ and $\Xi$, the following points ${\sigma}_{j+1}$ and ${\xi}_{j+1}$ are recursively determined in a way such that lasts
\begin{equation}
\max_{z\in\Sigma}|s_j(z)|=|s_j({\sigma}_{j+1})|\qquad\mbox{and}\qquad\min_{z\in\Xi}|s_j(z)|=|s_j({\xi}_{j+1})|.
\end{equation}

The points $\{({\sigma}_j,{\xi}_j)\}$ are called generalized Leja points.

We want to approximate the Markov function $z^{-1/2}$. In order to archive this, in \cite{guttel1} it is suggested to choose sets $\Sigma$ and $\Xi$ as $\Sigma=[{\lambda}_{min},{\lambda}_{max}]$ and $\Xi=[-\infty,0]$. These sets are optimal respectively for the choice of nodes and poles for that function.

\subsubsection{Adaptive Poles}
To approximate Markov functions, another way to choose the poles for the rational Arnoldi method is also qsuggested in \cite{guttel2}: they are called adaptive poles. We consider the function
\begin{equation}
s_m(z)=\frac{\prod_{k=1}^m(z-{\sigma}_k)}{q_{m-1}},
\end{equation}
in which ${\sigma}_k$ are the Ritz values of the projection of $A$ at the $m$-th step, i.e. the eigenvalues of $A_m$, and $q_{m-1}$ is the denominator associated with the corresponding rational Krylov space. Our aim is to minimize the error in the approximation for each step, and to do so we have to make $|s_m(z)|$ uniformly large on the set $\Xi$ choosing the next pole ${\xi}_m$ as
\begin{equation}
\min_{z\in\Xi}|s_m(z)|=|s_m({\xi}_m)|.
\end{equation}

This method is substantially a black-box method, because we do not need any information on the spectrum of $A$, but we only have to compute the eigenvalues of the matrix $A_m$, that are considerably less than the eigenvalues of $A$ itself.

\section{Experiments and Conclusions}
\label{experiments}
We are now going to consider some experiments in order to test the accuracy of the methods presented above. The first important thing is to analyse their convergence rate, that is how close they are to approximate the value of $(A\# B)v$. For this experiment we have used two positive definite matrices $A\in\mathbb{R}^{100\times 100}$ and $B\in\mathbb{R}^{100\times 100}$ and a vector $v\in\mathbb{R}^{100}$ randomly generated. We have run the methods up to 30 steps for each one.

The results are shown in the following graph, where we have on the x-axis the number of steps (or the number of poles considered for the Minimax method and the Gauss-Chebyshev quadrature, explained in \cite{trefethen} and \cite{iannazzo}, respectively) and on the y-axis the relative errors with respect to the real value of $(A\#B)v$ in a logarithmic scale.

\begin{figure}[htbp]
\centering
\includegraphics[scale=0.6]{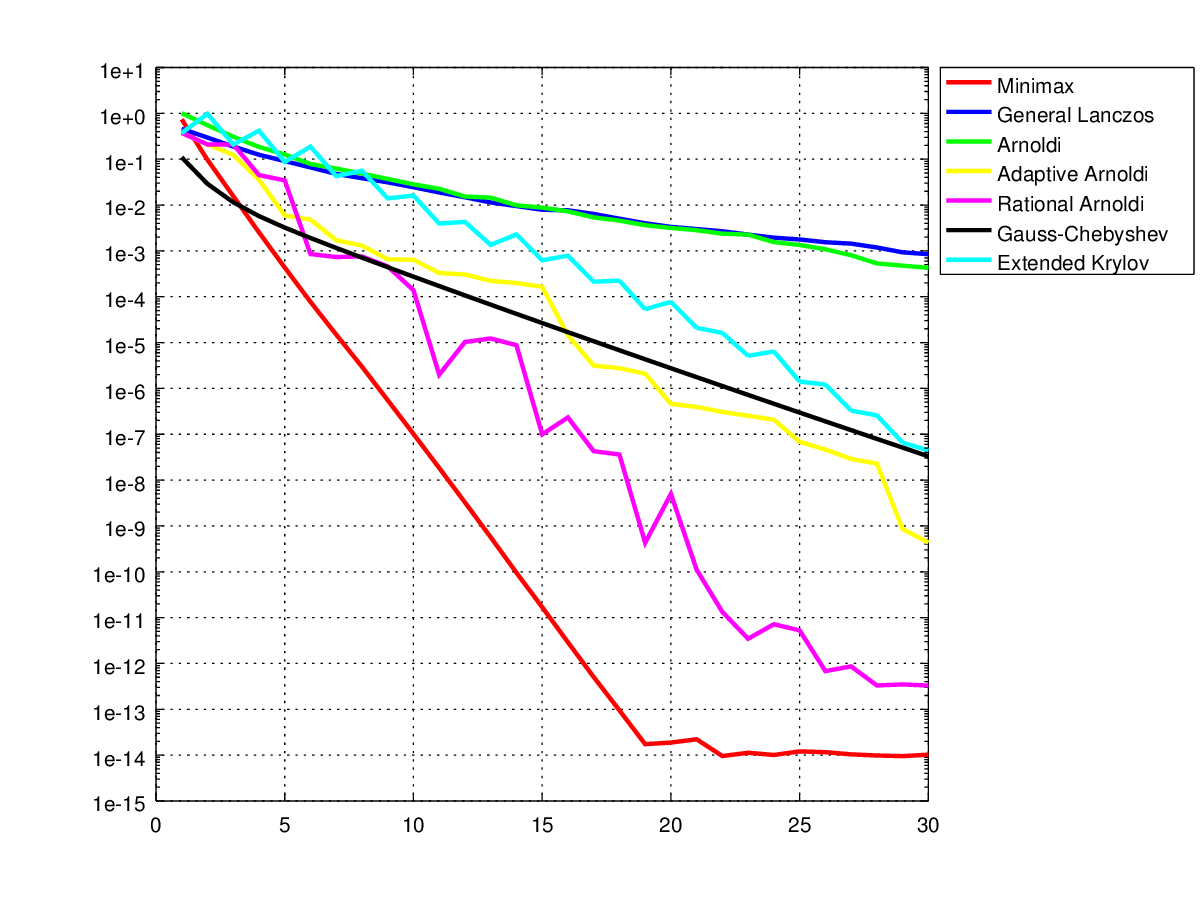}
\caption{Convergence rate for the presented methods.}
\label{fig:graphic}
\end{figure}

Figure~\ref{fig:graphic} shows how some methods, like the rational Arnoldi method with Leja points or adaptive poles, behave similarly to the Minimax method (they have almost the same convergence rate, with peaks in their graphics due to some numerical instabilities), but for the latter the knowledge of the full spectrum of the matrices is required.

Another important aspect we had to consider is how much time is required by these methods to compute the results. In fact, taking advantage of the sparse structure, Krylov methods use less operations than those methods that do not exploit it. We have used matrices $A$ and $B$ structured respectively as the finite-difference discretization matrices of the 1D and 2D Laplacian of variable dimension and a vector $v$ with every component equals to $1$. Both matrices are implemented using the Octave CSR sparse format.

$$A=\left[
\begin{array}{cccccc}
2 & -1 & & & & \\
-1 & 2 & -1 & & & \\
 & \ddots & \ddots & \ddots & & \\
 & & \ddots & \ddots & \ddots & \\
 & & & -1 & 2 & -1 \\
 & & & & -1 & 2
\end{array}
\right]
B=\left[
\begin{array}{cccccc}
4 & -1 & & -1 & & \\
-1 & 4 & -1 & & -1 & \\
 & \ddots & \ddots & \ddots & & -1 \\
-1 & & \ddots & \ddots & \ddots & \\
 & -1 & & -1 & 4 & -1 \\
 & & -1 & & -1 & 4
\end{array}
\right]$$

We have compared various methods: a method based on Schur decomposition presented in \cite{iannazzo} (that directly compute $A\#B$ and only later makes the product with $v$), the Minimax method presented in \cite{triplenick} (based on contour integrals) and two versions of the rational Arnoldi method, the one with Leja points and the adaptive one. Every method is run for 30 steps. We have measured the execution time of these methods with the Octave commands \texttt{tic} and \texttt{toc} (we have considered the time required to compute the extreme eigenvalues of $B^{-1}A$ with the Octave function \texttt{eigs()} in the estimation for the rational Arnoldi method with Leja points).

\begin{table}[htbp]
\begin{center}
\begin{tabular}{|c|c|c|c|c|}
\hline
Method & \multicolumn{4}{|c|}{Measured time} \\
\hline
Dimension & $1600\times 1600$ & $2500\times 2500$ & $3600\times 3600$ & $4900\times 4900$ \\
\hline
Minimax & 42.3826 & 136.398 & 382.048 & 743.067 \\
Schur & 24.3768 & 90.8797s & 253.466 & 651.758 \\
Leja Points & 2.26687 & 3.71883s & 6.18135 & 9.54854 \\
Adaptive & 2.03123 & 3.54691 & 5.5625 & 9.01572 \\
\hline
\end{tabular}
\end{center}
\caption{Time comparisons between the proposed methods.}
\label{fig:table}
\end{table}

Table~\ref{fig:table} shows how the computational time required by the Minimax method and the method based on Schur decomposition grows very quickly. For matrices of dimension $10000\times 10000$, the estimated times for the aforementioned methods become practically unsustainable, while methods based on Krylov spaces require only few seconds. These results show how important it is to preserve and exploit the sparse structure of the starting matrices in order to complete these computations in a reasonable time, gaining in terms of operations.

\bibliographystyle{plainnat}
\bibliography{Bibliografia}
\end{document}